\documentclass[10pt,a4paper]{article}
\usepackage[utf8]{inputenc}
\usepackage{amsmath}
\usepackage{amsfonts}
\usepackage{amssymb}
\usepackage{amsthm}
\usepackage{array}
\usepackage{hyperref}
\hypersetup{
    colorlinks,
    linkcolor={blue},
    citecolor={blue},
    urlcolor={blue}
}
\newtheorem{theorem}{Theorem}[section]
\newtheorem{lemma}[theorem]{Lemma}
\newtheorem{definition}[theorem]{Definition}
\newtheorem{example}[theorem]{Example}
\newtheorem{corollary}[theorem]{Corollary}
\usepackage{cleveref}
\usepackage[normalem]{ulem}
\author{Adrian Fellhauer}
\title{Approximation of smooth functions using Bernstein polynomials in multiple variables}
\begin{document}

\maketitle
\begin{abstract}
In this survey, we use (more or less) elementary means to establish the well-known result that for any given smooth multivariate function, the respective multivariate Bernstein polynomials converge to that function in all derivatives on each compact set. We then go on to strengthen that result to obtain that any smooth function on $\mathbb R^d$ may be approximated locally uniformly in all derivatives by \emph{one} sequence of polynomials. We will use neither the axiom of choice nor the power set axiom. We will use the method of proof by contradiction.
\end{abstract}

\tableofcontents

\section{Introduction}
The main question that concerns us in this survey is this: Given a function $f \in \mathcal C^\infty(\mathbb R^d)$, can we find a sequence of rational-coefficient polynomials $(p_n)_{n \in \mathbb N}$ in $d$ variables such that for each compact $K \subset \mathbb R^d$ and each multiindex $\beta \in \mathbb N_0^d$ (see \Cref{sec4}),
\begin{equation*}
\|\partial^\beta f - \partial^\beta p_n\|_{\mathcal C(K)} \overset{n \to \infty}{\longrightarrow} 0, ~~~~ \text{where } \| g \|_{\mathcal C(K)} := \sup_{x \in K} \|g(x)\| ~~?
\end{equation*}
As we will see, the answer is yes, and the main tool to help us are the \emph{Bernstein polynomials}, named after their inventor Sergei Natanovich Bernstein. For a given function $f \in \mathcal C([0, 1])$, they are defined as thus (see Bernstein \cite{bernstein}):
\begin{equation*}
B_n(f)(x) := \sum_{k=0}^n \binom{n}{k} f \left( \frac{n}{k} \right) x^k (1 - x)^{n-k}
\end{equation*}
We will generalize this definition to the multivariate case, prove that on $[0,1]^d$ the Bernstein polynomials converge to the respective function in all derivatives, then approximate the Bernstein polynomials by rational-coefficient polynomials and use the triangle inequality.

This survey is a successor of a draft, written by me for an educational website, which I abandoned (since it contained lots of flaws and was also otherwise of inferior quality).

\section{The classical Weierstraß approximation theorem and Bernstein's probabilistic proof}

A first step in the direction of affirmatively answering the question asked above is the Weierstrass approximation theorem, which states that any function $f \in \mathcal C([0, 1])$ can be arbitrarily approximated by polynomials in supremum norm. We shall state and prove this theorem following Bernstein \cite{bernstein}, except for fixing the slight flaw that only pointwise convergence was proven. The proof uses elementary probability theory.

\begin{theorem}
Let $f \in \mathcal C([0, 1])$. Then there exists a sequence $(p_n)_{n \in \mathbb N}$ of polynomials with real coefficients such that
\begin{equation*}
\lim_{n \to \infty} \|f - p_n\|_\infty = 0.
\end{equation*}
\end{theorem}

\begin{proof}
Consider the probability space $\Omega = [0,1]$ with the usual Borel measure as probability measure. For each $x \in [0,1]$, define the event
\begin{equation*}
A_x := [0,x].
\end{equation*}
For each $x \in [0,1]$, we define a couple of random variables on $\Omega^n$:
\begin{equation*}
X^x_k(\omega) := \begin{cases}
1 & \omega_k \in A_x \\
0 & \omega_k \notin A_x
\end{cases},~~~~k \in \{1, \ldots, n\},
\end{equation*}
\begin{equation*}
S^x_n(\omega) := \sum_{k=1}^n X^x_k(\omega),
\end{equation*}
\begin{equation*}
Y^x_n(\omega) := f \left( \frac{S^x_n(\omega)}{n} \right).
\end{equation*}
The expectation of $Y^x_n$ is given by
\begin{align*}
E(Y^x_n) & = \sum_{k=1}^n Y(k) P(S^x_n = k) \\
& = \sum_{k=1}^n f\left( \frac{k}{n} \right) \binom{n}{k} x^k (1 - x)^{n-k}.
\end{align*}
This is the famous Bernstein polynomial, and hence we see how it arose.

By assumption, $f$ is continuous, and by the Heine–Cantor theorem, $f$ even is uniformly continuous. Thus, for every $\epsilon > 0$, we can find a $\delta > 0$ such that
\begin{equation*}
\forall x, y \in [0,1] : |x - y| < 2\delta \Rightarrow |f(x) - f(y)| < \frac{\epsilon}{2}.
\end{equation*}
This fact of uniform continuity implies the following: Let $x \in [0,1]$ be arbitrary, and let $M_x$ and $m_x$ be the maximum resp. minimum of $f$ on $(x - \delta, x + \delta) \cap [0,1]$. Then
\begin{equation*}
\left( M_x - f(x) < \frac{\epsilon}{2} \right) \wedge \left( f(x) - m_x < \frac{\epsilon}{2} \right).
\end{equation*}
Let $x \in [0,1]$ be arbitrary. We compute the expectation of the random variable $\frac{S^x_n}{n}$:
\begin{equation*}
E\left( \frac{S^x_n}{n} \right) = \frac{1}{n} \sum_{k=1}^n E(X^x_k) = \frac{1}{n} n x = x.
\end{equation*}
Furthermore,
\begin{equation*}
E\left((X^x_k)^2\right) = \int_{\Omega^n} X^x_k(\omega)^2 d\omega = \int_{\Omega^n} X^x_k(\omega) d\omega = E\left(X^x_k\right)
\end{equation*}
and hence
\begin{align*}
\operatorname{Var}\left(\frac{S^x_n}{n}\right) & = E\left(\left(\frac{S^x_n}{n}\right)^2\right) - \left( E\left(\frac{S^x_n}{n}\right) \right)^2 \\
& = \frac{1}{n^2} \left( \sum_{j, k = 1 \atop k \neq j}^n E(X^x_j) E(X^x_k) + \sum_{k=1}^n E\left((X^x_k)^2\right) \right) - x^2 \\
& = \frac{1}{n} x (1-x),
\end{align*}
using the independence (and hence uncorrelatedness) of the random variables $X^x_1, \ldots, X^x_n$. Thus, by Markov's inequality,
\begin{equation*}
P \left( \left| \frac{S_n^x}{n} - x \right| \ge \delta \right) \le \frac{x (1-x)}{n \delta^2} \le \frac{1}{n \delta^2}.
\end{equation*}
Denoting (following Bernstein's original notation)
\begin{equation*}
\eta := P \left( \left| \frac{S_n^x}{n} - x \right| \ge \delta \right) ~~\text{and}~~ L := \|f\|_{\mathcal C([0,1])},
\end{equation*}
the definition of $Y^x_n$ and the monotonicity of expectation imply
\begin{equation*}
m_x(1 - \eta) - L \eta\le E(Y^x_n) \le M_x(1 - \eta) + L \eta.
\end{equation*}
Adding a zero on both sides yields
\begin{equation*}
f(x) + (m_x - f(x)) - \eta(m_x + L) \le E(Y^x_n) \le f(x) + (M_x - f(x)) - \eta(M_x - L),
\end{equation*}
and by inserting uniform continuity and our above result from Chebyshev's inequality, it becomes clear that we get a uniform approximation result ($x$ was arbitrary).
\end{proof}

Due to its simplicity, this proof (or rather, a simplified variant of it rid of the probabilistic formalism) is now taught in at least one course on approximation theory.

Most of the remainder of this survey is devoted to proving a stronger claim, namely the multivariate, derivatives result, without probability theory.

\section{Elementaries on polynomial interpolation}

Given a function $f: \mathbb R \to \mathbb R$ and $n$ distinct points $x_1, \ldots, x_n \in \mathbb R$, one can find a polynomial $\pi \in \mathbb R[x]$ of degree $n-1$ that \emph{interpolates} $f$ in the points $x_1, \ldots, x_n$, in the sense that $f(x_j) = \pi(x_j)$ for all $j \in \{1, \ldots, n\}$. This is done as thus.

\begin{definition}
Let $x_1, \ldots, x_n \in \mathbb R$ be pairwise distinct points. Then the \emph{Lagrange polynomials} with respect to these points are defined as
\begin{equation*}
L_j(x) := \prod_{k=1 \atop k \neq j}^n \frac{x - x_k}{x_j - x_k}, ~~~~j \in \{1, \ldots, n\}.
\end{equation*}
\end{definition}

With this definition, we can construct the (unique) interpolating polynomial of degree $n-1$ as thus.
\begin{theorem}
Let $x_1, \ldots, x_n \in \mathbb R$ be pairwise distinct points, and let $f: \mathbb R \to \mathbb R$ be given. Then the polynomial
\begin{equation*}
\pi(x) := \sum_{k=1}^n f(x_k) L_k(x)
\end{equation*}
is the unique polynomial of degree $n-1$ that interpolates $f$ at $x_1, \ldots, x_n$, that is, $f(x_j) = \pi(x_j)$, $j \in \{1, \ldots, n\}$. We call $\pi$ the \emph{interpolating polynomial} of $f$ at $x_1, \ldots, x_n$.
\end{theorem}

\begin{proof}
From the observation $L_j(x_m) = \delta_{j, m}$ for $j, m \in \{1, \ldots, n\}$ follows that $\pi$ is indeed an interpolating polynomial. It is of degree $n-1$ since it is the sum of the respective Lagrange polynomials, which are all of degree $n-1$. For uniqueness, suppose $\rho(x)$ is another polynomial of degree $n-1$ interpolating $f$ at $x_1, \ldots, x_n$. Then the polynomial $\pi - \rho$ also has degree $n-1$ and has $n$ distinct zeros. Therefore, by algebra, all its coefficients are zero.
\end{proof}

One could now conjecture that the interpolation polynomials are already enough to approximate a certain function, perhaps even in all derivatives. Unfortunately, in general, this is not the case, even when not considering any derivatives (although one can pick the points $x_1, \ldots, x_n$ such that approximation in supremum norm takes place as $n \to \infty$; we will not explicate that here). But when $f$ is sufficiently differentiable, we can at least get an expression for the error at a point to the right of all interpolation points (taken from Phillips \cite[p. 8, Theorem 1.1.2]{phillips}); an expression we will need in proving our approximation results.

\begin{theorem}\label{thm33}
Let $f: \mathbb R \to \mathbb R$ be $k+1$ times differentiable and let $x_1, \ldots, x_n \in \mathbb R$ be pairwise disjoint. Let further $y \in \mathbb R$ be another point strictly greater than any of the $x_1, \ldots, x_n$. Let $\pi \in \mathbb R[x]$ be the interpolating polynomial of $f$ at points $x_1, \ldots, x_n$. Then there exists a $\xi \in [x_1, y]$ such that
\begin{equation*}
f(y) - \pi(y) = \frac{\omega(y)f^{(n)}(\xi)}{n!},
\end{equation*}
where
\begin{equation*}
\omega(x) := (x - x_1) \cdots (x - x_n)
\end{equation*}
is the monic polynomial with roots $x_1, \ldots, x_n$
\end{theorem}

\begin{proof}
We define the function
\begin{equation*}
G(t) := f(t) - \pi(t) - \frac{\omega(t) \left(f(y) - \pi(y)\right)}{\omega(y)}.
\end{equation*}
This function has $n+1$ zeroes, namely at $x_1, \ldots, x_n, y$. Hence, by repeated application of Rolle's theorem, $G^{(n)}$ has a root which we shall denote by $\xi$. But on the other hand, taking into account the form of the polynomials $\pi$ and $\omega$, $G^{(n)}$ evaluates to
\begin{equation*}
G^{(n)}(t) = f^{(n)}(t) - \frac{n!\left(f(y) - \pi(y)\right)}{\omega(y)},
\end{equation*}
and hence reformulating the equation $G^{(n)}(\xi) = 0$ gives the theorem.
\end{proof}

Below, this error expression will allow us to derive an alternative formula for the derivatives of the Bernstein polynomials, which will ultimately allow us to prove the approximative property we aim at.

\section{Multiindex notation}\label{sec4}

In this section, we establish several notational conventions for multiindices in order to be able to briefly write down some required formulae.\footnote{These notations were largely invented by Laurent Schwartz \cite[p.14f.]{schwartz}.} Recall that a multiindex of dimension $d \in \mathbb N$ is a vector $(\alpha_1, \ldots, \alpha_d) \in \mathbb N_0^d$. We may write partial derivatives as thus:
\begin{equation*}
\frac{\partial^{\alpha_1}}{\partial x_1} \cdots \frac{\partial^{\alpha_d}}{\partial x_d} f =: \partial^\alpha f;
\end{equation*}
according to Clairaut's theorem, by expressions like this we may express all possible partial derivatives of any order. We have a multiindex binomial coefficient as thus:
\begin{equation*}
\binom{\alpha}{\beta} := \binom{\alpha_1}{\beta_1} \cdots \binom{\alpha_d}{\beta_d}.
\end{equation*}
We may define a partial order on multiindices:
\begin{equation*}
\alpha \le \beta :\Leftrightarrow \forall j \in \{1, \ldots, d\}: \alpha_j \le \beta_j.
\end{equation*}
If $n \in \mathbb N$ is any natural number, we denote the vector that is constantly $n$ in boldface:
\begin{equation*}
\mathbf n := (n, \ldots, n);
\end{equation*}
for instance, $\mathbf 1 = (1, \ldots, 1)$. If $\alpha \ge \mathbf 1$, we may define
\begin{equation*}
\frac{\beta}{\alpha} := \left( \frac{\beta_1}{\alpha_1}, \ldots, \frac{\beta_d}{\alpha_d} \right) \in \mathbb Q^d.
\end{equation*}
The factorial of a multiindex is
\begin{equation*}
\alpha! := \alpha_1! \cdots \alpha_d!.
\end{equation*}
The absolute value of a multiindex $\alpha$ is defined as
\begin{equation*}
|\alpha| := \sum_{k=1}^n \alpha_k.
\end{equation*}
The minimum of two multiindices $\alpha, \beta$ is
\begin{equation*}
\min\{\alpha, \beta\} := \left( \min\{\alpha_1, \beta_1\}, \ldots, \min\{\alpha_d, \beta_d\} \right)
\end{equation*}
Further, we write $f \in \mathcal C^\beta(O)$, $O \subseteq \mathbb R^d$ open, if on $O$ all the derivatives
\begin{equation*}
\partial^\beta f
\end{equation*}
exist and are continuous.

\section{Approximation in one dimension}

In this section, we shall prove that any $f \in \mathcal C^\infty([0, 1])$ can be approximated in all derivatives by one sequence of polynomials; that is, there is a sequence $(p_n)_{n \in \mathbb N}$ of polynomials such that $\|\partial^\alpha f - \partial^\alpha p_n\|_{\mathcal C([0,1])} \to 0$ for all multiindices $\alpha$. We do this also because later, we can elevate this result to the multivariate case, and then we will strengthen the result in the sense that we will prove that we may pick the approximating sequence of polynomials as \emph{rational} polynomials and even get locally uniform convergence on $\mathbb R^d$. We will mostly follow the treatise of Phillips \cite[Chapter 7, p. 247 - 258]{phillips}, except for not using forward differences and justifying a certain sum manipulation in a different, perhaps more intuitive and universally applicable form. For the theorem, we need a couple of fairly elementary and technical lemmata, which are exercises in elementary analysis. First, we shall prove an elementary identity which follows from the Binomial theorem and sum manipulation.

\begin{lemma}
For all $N \in \mathbb N$, $N \le n$, we have
\begin{equation*}
\sum_{k=0}^n k^{\underline{N}} \binom{n}{k} x^k (1 - x)^{n-k} = n^{\underline{N}} x^N.\footnote{Note that $a^{\underline{b}} = a (a-1) \cdots (a - b + 1)$, $a, b \in \mathbb N_0$ is the falling factorial.}
\end{equation*}
\end{lemma}

\begin{proof}
We first note that
\begin{equation*}
\sum_{k=0}^n k^{\underline{N}} \binom{n}{k} x^k (1 - x)^{n-k} = \sum_{k=N}^n k^{\underline{N}} \binom{n}{k} x^k (1 - x)^{n-k}
\end{equation*}
since the first $N$ summands just vanish by definition of the falling factorial. We further compute
\begin{align*}
\sum_{k=N}^n k^{\underline{N}} \binom{n}{k} x^k (1 - x)^{n-k} & = \sum_{k=N}^n \frac{k!}{(k-N)!} \frac{n!}{k!(n-k)!} x^k (1-x)^{n-k} \\
& = n^{\underline{N}} x^N \sum_{k=N}^n \frac{(n-N)!}{(k-N)!((n-N)-(k-N))!} x^{k-N} (1-x)^{(n-N)-(k-N)} \\
& = n^{\underline{N}} x^N \sum_{k=0}^{n-N} \binom{n-N}{k} x^k (1-x)^{(n-N)-k},
\end{align*}
and since by the binomial theorem
\begin{equation*}
1 = (x + (1 - x))^{N-n} = \sum_{k=0}^{n-N} \binom{n-N}{k} x^k (1-x)^{(n-N)-k},
\end{equation*}
the claim follows.
\end{proof}

When plugging in different $N$, the above lemma spits out several formulae. In fact, for our purposes we will only need the cases $N = 0, 1, 2$, which look like this:

\begin{example}\label{ex52}\[\def\arraystretch{3.0}\begin{array}[t]{@{~~~~~~} >{\displaystyle}r @{~~~~~~~~} >{\displaystyle}r @{~} >{\displaystyle}l}
N = 0: & \sum_{k=0}^n \binom{n}{k} x^k (1 - x)^{n-k} &= 1 \\
N = 1: & \sum_{k=0}^n k \binom{n}{k} x^k (1 - x)^{n-k} &= n x \\
N = 2: & \sum_{k=0}^n k (k-1) \binom{n}{k} x^k (1 - x)^{n-k} &= n(n-1)x^2 \\
\end{array}
\]

~\,
\end{example}

\begin{corollary}
\begin{equation*}
\sum_{k=0}^n (k - nx)^2 \binom{n}{k} x^k (1 - x)^{n-k} = n x (1-x).
\end{equation*}
\end{corollary}

\begin{proof}
By the formulae in \Cref{ex52},
\begin{align*}
\sum_{k=0}^n (k - nx)^2 \binom{n}{k} x^k (1 - x)^{n-k} & = \sum_{k=0}^n (k^2 + n^2 x^2 - 2knx) \binom{n}{k} x^k (1 - x)^{n-k} \\
& = \sum_{k=0}^n (k(k-1) + k + n^2 x^2 - 2knx) \binom{n}{k} x^k (1 - x)^{n-k} \\
& = n(n-1)x^2 + nx + n^2 x^2 - 2n^2x^2 = n x (1-x).
\end{align*}
\end{proof}

We also derive a certain expression for the derivatives of Bernstein polynomials (of sufficient order), which will be crucial in our proof of convergence.

\begin{lemma}\label{lem54}
Let $f \in \mathcal C^j([0,1])$ and let $n \in \mathbb N$. Then there exist $\xi_k \in \left[ \frac{k}{n+j}, \frac{k+j}{n+j} \right]$ such that
\begin{equation*}
\frac{d^j}{dx^j} B_{n+j}(f)(x) = \frac{(n+j)!}{n!(n+j)^j} \sum_{k=1}^n \binom{n}{k} f^{(j)}(\xi_k) x^k (1 - x)^{n-k}.
\end{equation*}
\end{lemma}

\begin{proof}
By the general Leibniz rule (see, for instance, Forster \cite[Exercise 15.11 i), p. 168]{forster}),
\scriptsize
\begin{align*}
\frac{d^j}{dx^j} B_{n+j}(f)(x) & = \frac{d^j}{dx^j} \sum_{k=0}^{n+j} \binom{n+j}{k} f\left( \frac{k}{n+j} \right) x^k (1-x)^{n+j-k} \\
& = \sum_{k=0}^{n+j} \binom{n+j}{k} f\left( \frac{k}{n+j} \right) \frac{d^j}{dx^j} x^k (1-x)^{n+j-k} \\
& = \sum_{k=0}^{n+j} \binom{n+j}{k} f\left( \frac{k}{n+j} \right) \sum_{m=\max \{0, k-n\}}^{\min\{k,j\}} \binom{j}{m} k^{\underline{m}} x^{k-m} (-1)^{j-m} (n+j-k)^{\underline{j-m}} (1-x)^{n+j-k-(j-m)} \\
& = \sum_{k=0}^{n+j} \sum_{m=\max \{0, k-n\}}^{\min\{k,j\}} f\left( \frac{k}{n+j} \right) \binom{n+j}{k} \binom{j}{m} k^{\underline{m}} (n+j-k)^{\underline{j-m}} x^{k-m} (-1)^{j-m} (1-x)^{n+j-k-(j-m)}.
\end{align*}\normalsize
Further,
\begin{align*}
\binom{n+j}{k} k^{\underline{m}} (n+j-k)^{\underline{j-m}} & = \frac{(n+j)!}{k!(n+j-k)!} \frac{k!}{(k-m)!} \frac{(n+j-k)!}{(n+j-k-(j-m))!} \\
& = \frac{(n+j)!}{n!} \frac{n!}{(k-m)!(n+j-k-(j-m))!} \\
& = \frac{(n+j)!}{n!} \binom{n}{k-m},
\end{align*}
which is why\small
\begin{equation*}
\frac{d^j}{dx^j} B_{n+j}(f)(x) = \frac{(n+j)!}{n!} \sum_{k=0}^{n+j} \sum_{m=\max \{0, k-n\}}^{\min\{k,j\}} \binom{n}{k-m} f\left( \frac{k}{n+j} \right) \binom{j}{m} x^{k-m} (-1)^{j-m} (1-x)^{n+j-k-(j-m)}.
\end{equation*}\normalsize
To proceed further, we need to notice the following summation technique justified by elementary set theory. Assume we have two sets $S_1$ and $S_2$ and two functions $f_j: S_j \to \mathbb R, j = 1, 2$. We consider the sums
\begin{equation*}
\sum_{s \in S_1} f_1(s) ~~\text{and}~~ \sum_{t \in S_2} f_2(t).
\end{equation*}
Assume that there is a bijection $\Phi: S_1 \to S_2$ such that $\forall s \in S_1: f_1(s) = f_2(\Phi(s))$. Then the two sums are equal, because one sums the same numbers.

Applying this argument for
\begin{align*}
S_1 &= \left\{ (k,m) \middle| \left( 0 \le k \le n+j \right) \wedge \left( \max \{0, k-n\} \le m \le \min\{k,j\} \right) \right\}, \\
S_2 &= \left\{ (k,m) \middle| \left( 0 \le k \le n \right) \wedge \left( 0 \le m \le j \right) \right\},\\
f_1(k,m) & = \binom{n}{k-m} f\left( \frac{k}{n+j} \right) \binom{j}{m} x^{k-m} (-1)^{j-m} (1-x)^{n+j-k-(j-m)},\\
f_2(k,m) & = \binom{n}{k} (-1)^{j-m} f\left( \frac{k+m}{n+j} \right) \binom{j}{m} x^k (1-x)^{n-k}
\end{align*}
and
\begin{equation*}
\Phi(k,m) = (k-m,m)
\end{equation*}
we get
\begin{equation*}
\frac{d^j}{dx^j} B_{n+j}(f)(x) = \frac{(n+j)!}{n!} \sum_{k=0}^n \binom{n}{k} \sum_{m=0}^j (-1)^{j-m} f\left( \frac{k+m}{n+j} \right) \binom{j}{m} x^k (1-x)^{n-k}.
\end{equation*}
Now the connection to interpolation polynomials comes in; for, we claim
\begin{equation*}
\sum_{m=0}^j (-1)^{j-m} f\left( \frac{k+m}{n+j} \right) \binom{j}{m} = f \left( \frac{k+j}{n+j} \right) - \pi \left( \frac{k+j}{n+j} \right),
\end{equation*}
where $\pi \in \mathbb R[x]$ interpolates $f$ at the points $\frac{k}{n+j}, \ldots, \frac{k+j-1}{n+j}$. Indeed, as shown in chapter~3, $\pi$ is a linear combination of the Lagrange polynomials at these points as thus:
\begin{align*}
\pi(x) & = \sum_{m=1}^j f \left( \frac{k + m - 1}{n+j} \right) L_m (x) \\
& = \sum_{m=0}^{j-1} f \left( \frac{k + m}{n+j} \right) L_{m+1} (x),
\end{align*}
where the Lagrange polynomials $L_m$ are defined according to the points $\frac{k}{n+j}, \ldots, \frac{k+j-1}{n+j}$. But
\begin{align*}
L_{m+1} \left( \frac{k+j}{n+j} \right) & = \prod_{l=0 \atop l \neq m}^{j-1} \frac{\frac{k+j}{n+j} - \frac{k+l}{n+j}}{\frac{k+m}{n+j} - \frac{k+l}{n+j}} \\
& = \prod_{l=0 \atop l \neq m}^{j-1} \frac{j-l}{m-l} \\
& = \frac{j!}{j-m} \prod_{l=0}^{m-1} \frac{1}{m-l} \prod_{l=m+1}^{j-1} \frac{1}{m-l} \\
& = j! \frac{(-1)^{j-m+1}}{m!(j-m)!}
\end{align*}
which is why the claim is really true. Now by \Cref{thm33},
\begin{equation*}
f \left( \frac{k+j}{n+j} \right) - \pi \left( \frac{k+j}{n+j} \right) = \frac{f^{j}(\xi_k) \omega\left( \frac{k+j}{n+j} \right)}{j!}
\end{equation*}
for a certain $\xi_k \in \left[ \frac{k}{n+j}, \frac{k+j}{n+j} \right]$ and $\omega$ defined as in \Cref{thm33}. But
\begin{align*}
\omega\left( \frac{k+j}{n+j} \right) & = \prod_{l=1}^j \left( \frac{k+j}{n+j} - \frac{k+j-l}{n+j} \right) \\
& = \frac{j!}{(n+j)^j}
\end{align*}
and hence
\begin{equation*}
\frac{d^j}{dx^j} B_{n+j}(f)(x) = \frac{(n+j)!}{n!(n+j)^j} \sum_{k=1}^n \binom{n}{k} f^{(j)}(\xi_k) x^k (1 - x)^{n-k}
\end{equation*}
as desired.
\end{proof}

\begin{corollary}\label{cor55}
Let $j \in \mathbb N_0$, $f \in \mathcal C^j([0, 1])$. Then the family
\begin{equation*}
\left( \frac{d^j}{dx^j} B_{n+j}(f)(x) \right)_{n \in \mathbb N}
\end{equation*}
is uniformly bounded on $[0, 1]$; in fact,
\begin{equation*}
\forall n \in \mathbb N: \left\| \frac{d^j}{dx^j} B_{n+j}(f) \right\|_{\mathcal C([0, 1])} \le \|f^{(j)}\|_{\mathcal C([0, 1])}.
\end{equation*}
\end{corollary}

\begin{proof}
By the previous lemma and in the last line using the case $N=0$ from \Cref{ex52} we get
\begin{align*}
\left\| \frac{d^j}{dx^j} B_{n+j}(f) \right\|_{\mathcal C([0, 1])} & = \sup_{x \in [0, 1]} \left| \frac{(n+j)!}{n!(n+j)^j} \sum_{k=1}^n \binom{n}{k} f^{(j)}(\xi_k) x^k (1 - x)^{n-k} \right| \\
& \le \sup_{x \in [0, 1]} \underbrace{\frac{(n+j)!}{n!(n+j)^j}}_{\le 1} \sum_{k=1}^n |f^{(j)}(\xi_k)| \binom{n}{k} x^k (1 - x)^{n-k} \\
& \le \|f^{(j)}\|_{\mathcal C([0, 1])};
\end{align*}
note that the $\xi_k$ depend on $x$ but it doesn't matter since we only used the estimate $|f^{(j)}(\xi_k)| \le \|f^{(j)}\|_{\mathcal C([0, 1])}$ which doesn't depend on the $\xi_k$.
\end{proof}

Now we're ready to formulate and prove the one-dimensional approximation theorem, which shall later be the cornerstone of the proof of the main result.

\begin{theorem}\label{thm56}
Let $f \in \mathcal C^j([0, 1])$. Then
\begin{equation*}
\left\|\frac{d^j}{dx^j} f - \frac{d^j}{dx^j} B_n(f)\right\|_{\mathcal C([0, 1])} \to 0, n \to \infty.
\end{equation*}
\end{theorem}

\begin{proof}
We prove instead
\begin{equation*}
\left\|\frac{d^j}{dx^j} f - \frac{d^j}{dx^j} B_{n+j}(f)\right\|_{\mathcal C([0, 1])} \to 0, n \to \infty.
\end{equation*}
By the triangle inequality,
\begin{align*}
\left\|\frac{d^j}{dx^j} f - \frac{d^j}{dx^j} B_{n+j}(f)\right\|_{\mathcal C([0, 1])} \le & \left\|\frac{d^j}{dx^j} f - \frac{n!(n+j)^j}{(n+j)!} \frac{d^j}{dx^j} B_{n+k}(f)\right\|_{\mathcal C([0, 1])}\\
& + \left| 1 - \frac{n!(n+j)^j}{(n+j)!} \right| \left\|\frac{d^j}{dx^j} B_{n+k}(f)\right\|_{\mathcal C([0, 1])}.
\end{align*}
The latter summand on the right converges to zero by \Cref{cor55}. So all that is left to show is that the first summand on the right converges to zero as well. Let thus $\epsilon > 0$. We want to prove that if $n \ge N$ for a sufficiently large $N$,
\begin{equation*}
\left| \frac{d^j}{dx^j} f(x) - \frac{n!(n+j)^j}{(n+j)!} \frac{d^j}{dx^j} B_{n+k}(f)(x) \right| < \epsilon
\end{equation*}
for all of $x \in [0, 1]$. First, we manipulate the expression inside the absolute value. Indeed, using $N=0$ in \Cref{ex52} and \Cref{lem54} we get
\begin{equation*}
\frac{d^j}{dx^j} f(x) - \frac{n!(n+j)^j}{(n+j)!} \frac{d^j}{dx^j} B_{n+k}(f)(x) = \sum_{k=0}^n \binom{n}{k} \left( f^{(j)}(x) - f^{(j)}(\xi_k) \right) x^k (1-x)^{n-k}.
\end{equation*}
for suitable $\xi_k \in \left[ \frac{k}{n +j}, \frac{k + j}{n +j} \right]$. Thus, by the triangle inequality,
\begin{equation*}
\left| \frac{d^j}{dx^j} f(x) - \frac{n!(n+j)^j}{(n+j)!} \frac{d^j}{dx^j} B_{n+k}(f)(x) \right| \le \sum_{k=0}^n \binom{n}{k} \left| f^{(j)}(x) - f^{(j)}(\xi_k) \right| x^k (1-x)^{n-k}.
\end{equation*}
By assumption, $f^{(j)}$ is continuous, and even uniformly continuous on the compact set $[0,1]$ by the Heine–Cantor theorem. Thus, pick $\delta > 0$ such that
\begin{equation*}
\forall y, z \in [0,1]: |y - z| < 2 \delta \Rightarrow \left|f^{(j)}(y) - f^{(j)}(z)\right| < \frac{\epsilon}{2}.
\end{equation*}
We have
\begin{align*}
\sum_{k=0}^n \binom{n}{k} \left| f^{(j)}(x) - f^{(j)}(\xi_k) \right| x^k (1-x)^{n-k} = & \sum_{k=0 \atop \left| x - \frac{k}{n+j} \right| < \delta}^n \binom{n}{k} \left| f^{(j)}(x) - f^{(j)}(\xi_k) \right| x^k (1-x)^{n-k} \\
& + \sum_{k=0 \atop \left| x - \frac{k}{n+j} \right| \ge \delta}^n \binom{n}{k} \left| f^{(j)}(x) - f^{(j)}(\xi_k) \right| x^k (1-x)^{n-k}
\end{align*}
and for sufficiently large $n$ (say $n \ge N_1$), the first sum is bounded as in
\begin{equation*}
\sum_{k=0 \atop \left| x - \frac{k}{n+j} \right| < \delta}^n \binom{n}{k} \left| f^{(j)}(x) - f^{(j)}(\xi_k) \right| x^k (1-x)^{n-k} \le \frac{\epsilon}{2}
\end{equation*}
by means of uniform continuity and $N=0$ in \Cref{ex52}. Now if $\left| x - \frac{k}{n+j} \right| \ge \delta$, then
\begin{equation*}
\frac{1}{(n+j)\delta} |(n+j)x - k| \ge 1 \text{ and } \frac{1}{(n+j)^2 \delta^2}\left( (n+j)x - k \right)^2 \ge 1
\end{equation*}
and therefore by \Cref{cor55},
\begin{align*}
\sum_{k=0 \atop \left| x - \frac{k}{n+j} \right| \ge \delta}^n \binom{n}{k} \left| f^{(j)}(x) - f^{(j)}(\xi_k) \right| x^k (1-x)^{n-k} & \le \frac{2 \|f^{(j)}\|_\infty}{(n+j)^2 \delta^2} \sum_{k=0}^n \left( (n+j)x - k \right)^2 \binom{n}{k} x^k (1-x)^{n-k} \\
& = \frac{2 \|f^{(j)}\|_\infty}{(n+j)^2 \delta^2} (n+j)x(1-x).
\end{align*}
If we pick $N_2$ large enough such that the last expression is less than $\frac{\epsilon}{2}$ for $n \ge N_2$, and then set $N = \max\{N_1, N_2\}$, we have completed our proof, since collecting things together, we have found that the sum we wanted to bound by $\epsilon$ is in fact bound by $\epsilon$.
\end{proof}

\section{Approximation in multiple dimensions}

\begin{definition}
If $\alpha \ge \mathbf 1$, by analogy we define multivariate Bernstein polynomials as thus (see e.g. \cite[(5.11), p.119]{reimer}):
\begin{equation*}
B_\alpha(f)(x) := \sum_{\beta \in \mathbb N_0^d \atop \beta \le \alpha} \binom{\alpha}{\beta} f\left( \frac{\beta}{\alpha} \right) x^\beta (\mathbf 1 - x)^{\alpha - \beta}.
\end{equation*}
\end{definition}

We will also define multivariate Bernstein polynomials for $\alpha \not\ge \mathbf 1$, but we need to be a bit careful because in these cases, we want to achieve that if $\alpha_k$ is zero, then $f$ shall not be ``Bernstein-expanded'' in direction $k$. Hence, the following definition is suitable:

\begin{definition}
Let $\alpha \in \mathbb N_0^d$, $\alpha \neq 0$ and let $k_1, \ldots, k_m \in \{1, \ldots,d\}$ be the indices where $\alpha$ is nonzero. Then we define
\begin{equation*}
B_\alpha(f)(x) := \sum_{\beta \in \mathbb N_0^m \atop \beta \le \alpha^+} \binom{\alpha^+}{\beta} f \left( x_1, \ldots, \frac{\beta_1}{\alpha_{k_1}}, x_{k_1+1}, \ldots, \frac{\beta_2}{\alpha_{k_2}}, \ldots, \frac{\beta_m}{\alpha_{k_m}}, \ldots, x_n \right) x^{\beta_0} (\mathbf 1 - x)^{\alpha - \beta_0}
\end{equation*}
where $\alpha^+ := (\alpha_{k_1}, \ldots, \alpha_{k_m})$, $\beta_0 := \left(0, \ldots, \beta_1, 0, \ldots, \beta_2, \ldots, \beta_m, \ldots, 0 \right)$, where $\beta_j$ is at place $k_j$.
\end{definition}

In fact, we will only need this definition in the cases where only the first component is nonzero and where only the first component is zero. This already suffices to get the induction going, since we will split $\alpha$ up in a suitable way. The following lemma makes this more explicit.
\begin{lemma}\label{lemma63}
Let $\alpha \ge \mathbf 1$, $\beta \in \mathbb N_0^d$, $\alpha = (\alpha_1, \ldots, \alpha_d)$, $\beta = (\beta_1, \ldots, \beta_d)$ and $f \in \mathcal C^\beta(\mathbb R^d)$. Set $\alpha' := (0, \alpha_2, \ldots, \alpha_d)$, $\beta' := (0, \beta_2, \ldots, \beta_d)$. Then
\begin{equation*}
\partial^\beta B_\alpha(f) = \partial^{\beta'} B_{\alpha'} \left( \partial^{\beta_1 \cdot e_1} B_{\alpha_1 \cdot e_1} (f) \right),
\end{equation*}
where as usual $e_1 = (1, 0, \ldots, 0)$ denotes the first vector of the standard basis of $\mathbb R^d$.
\end{lemma}

\begin{proof}
\small
\begin{align*}
\partial^{\beta'} B_{\alpha'} \left( \partial^{\beta_1 \cdot e_1} B_{\alpha_1 \cdot e_1} (f) \right)(x) & = \partial^{\beta'} \sum_{\gamma \in \mathbb N_0^{d-1} \atop \gamma \le \alpha^+} \partial^{\beta_1 \cdot e_1} B_{\alpha_1 \cdot e_1} (f) \left(x_1, \frac{\gamma_1}{\alpha_2}, \ldots, \frac{\gamma_{d-1}}{\alpha_d} \right) \binom{\alpha^+}{\gamma} x^{\gamma_0} (\mathbf 1 - x)^{\alpha' - \gamma_0} \\
& = \partial^\beta \sum_{\gamma \in \mathbb N_0^{d-1} \atop \gamma \le \alpha^+} \sum_{k=0}^{\alpha_1} \binom{\alpha^+}{\gamma} \binom{\alpha_1}{k} f \left( \frac{k}{\alpha_1}, \frac{\gamma_1}{\alpha_2}, \ldots, \frac{\gamma_{d-1}}{\alpha_d} \right) x^{\gamma_0} (\mathbf 1 - x)^{\alpha' - \gamma_0} x_1^k (1-x_1)^{\alpha_1 - k} \\
& = \partial^\beta B_\alpha(f),
\end{align*}\normalsize
where $\alpha^+ =(\alpha_2, \ldots, \alpha_d)$.
\end{proof}
Similar identities hold true as well, and will be used in proving \Cref{thm64}. To derive them is equally challenging.

The (multivariate) Bernstein polynomials have several nice properties which will later help us with the proof of \Cref{thm67}. First of all, all their derivatives are bound by the respective derivatives of the respective functions:

\begin{theorem}\label{thm64}
Let $f \in \mathcal C^\beta(\mathbb R^d)$, $\beta \in \mathbb N_0^d$. Then
\begin{equation*}
\forall \alpha > \beta : \|\partial^\beta B_\alpha (f)\|_{\mathcal C([0,1]^d)} \le \|\partial^\beta f\|_{\mathcal C([0,1]^d)}.
\end{equation*}
\end{theorem}

Here, $\alpha > \beta$ means $\alpha_k > \beta_k$ for all $k$.

\begin{proof}
We prove the theorem by induction on $d$. The case $d=1$ is given by \Cref{cor55}. Let the case $d-1$ be proven. First,  for fixed $x_1 \in [0,1]$ define $\mu_{x_1}(x_2, \ldots, x_d) := \partial^{\beta_1 \cdot e_1} B_{\alpha_1 \cdot e_1} (f)(x_1, x_2, \ldots, x_d)$ and for fixed $(x_2, \ldots, x_d) \in [0,1]^{d-1}$ define $\lambda_{(x_2, \ldots, x_d)}(x_1) := \partial^{\beta'} f(x_1, x_2, \ldots, x_d)$. Note the identities
\begin{align*}
\partial^{\beta'} B_{\alpha'} \left( \partial^{\beta_1 \cdot e_1} B_{\alpha_1 \cdot e_1} (f) \right)(x_1, x_2, \ldots, x_d) & = \partial^{(\beta')^+} B_{(\alpha')^+} \left( \mu_{x_1} \right)(x_2, \ldots, x_d) ~~\text{and} \\
\partial^{(\beta')^+} \mu_{x_1}(x_2, \ldots, x_d) & = \partial^{\beta_1 \cdot e_1} B_{\alpha_1 \cdot e_1} (\lambda_{(x_2, \ldots, x_d)})(x_1).
\end{align*}
By the inductive hypothesis, for all $x_1 \in [0,1]$
\begin{equation*}
\left\| \partial^{(\beta')^+} B_{(\alpha')^+} \left( \mu_{x_1} \right) \right\|_{\mathcal C([0,1]^{d-1})} \le \left\| \partial^{(\beta')^+} \mu_{x_1} \right\|_{\mathcal C([0,1]^{d-1})},
\end{equation*}
and by the case $d=1$, for all $(x_2, \ldots, x_d) \in [0,1]^{d-1}$
\begin{equation*}
\left\| \partial^{\beta_1 \cdot e_1} B_{\alpha_1 \cdot e_1} (\lambda_{(x_2, \ldots, x_d)}) \right\|_{\mathcal C([0,1])} \le \left\| \partial^{\beta_1 \cdot e_1} \lambda_{(x_2, \ldots, x_d)} \right\|_{\mathcal C([0,1])}.
\end{equation*}
By \Cref{lemma63},
\begin{equation*}
\partial^\beta B_\alpha(f) = \partial^{\beta'} B_{\alpha'} \left( \partial^{\beta_1 \cdot e_1} B_{\alpha_1 \cdot e_1} (f) \right).
\end{equation*}
Hence,
\begin{align*}
\left\| \partial^\beta B_\alpha(f) \right\|_{\mathcal C([0,1]^d)} & = \sup_{x_1 \in [0,1]} \left\| \partial^{(\beta')^+} B_{(\alpha')^+} \left( \mu_{x_1} \right) \right\|_{\mathcal C([0,1]^{d-1})} \le \sup_{x_1 \in [0,1]} \left\| \partial^{(\beta')^+} \mu_{x_1} \right\|_{\mathcal C([0,1]^d)} \\
& = \sup_{(x_2, \ldots, x_d) \in [0,1]^{d-1}} \left\| \partial^{\beta_1 \cdot e_1} B_{\alpha_1 \cdot e_1} (\lambda_{(x_2, \ldots, x_d)}) \right\|_{\mathcal C([0,1])} \\
& \le \sup_{(x_2, \ldots, x_d) \in [0,1]^{d-1}} \left\| \partial^{\beta_1 \cdot e_1} \lambda_{(x_2, \ldots, x_d)} \right\|_{\mathcal C([0,1])} \\
& = \|\partial^\beta f\|_{\mathcal C([0,1]^d)}.
\end{align*}
\end{proof}

Roughly speaking, this theorem tells us that the derivatives of the Bernstein polynomials with respect to a function are not more ``rough'' than the derivatives of the function itself.

As a corollary, we get:

\begin{corollary}\label{cor65}
Let $f \in \mathcal C^\beta(\mathbb R^d)$, $\beta \in \mathbb N_0^d$. Then the family of functions
\begin{equation*}
\left( \partial^\beta B_\alpha (f) \right)_{\alpha > \beta}
\end{equation*}
is equicontinuous.
\end{corollary}

\begin{theorem}
\begin{equation*}
B_\alpha(f + \lambda g) = B_\alpha(f) + \lambda B_\alpha(g).
\end{equation*}
\end{theorem}

Hence, \Cref{thm64} can be reinterpreted as thus: The linear operator 
\begin{equation*}
\mathcal C^\beta([0,1]^d) \to \mathcal C^\beta([0,1]^d), f \mapsto \partial^\beta B_\alpha(f)
\end{equation*}
is continuous with Lipschitz constant $1$.

\begin{theorem}\label{thm67}
Let $f \in \mathcal C^\beta(\mathbb R^d)$. Then
\begin{equation*}
\|\partial^\beta B_\alpha(f) - \partial^\beta f\|_{\mathcal C([0,1]^d)} \to 0 ~~\text{as}~~ \min_{j \in \{1, \ldots, d\}} \alpha_j \to \infty.
\end{equation*}
\end{theorem}

\begin{proof}
We prove the theorem by induction on $d$. $d = 1$ is \Cref{thm56}. Assume the case $d-1$ be given. We start off by proving the following slightly different claim:

\emph{Let $f \in \mathcal C^\beta(\mathbb R^d)$. If $\beta$ is a multiindex such that \emph{at least} one of the entries of $\beta$ is zero, then
\begin{equation*}
\|\partial^\beta B_\alpha(f) - \partial^\beta f\|_{\mathcal C([0,1]^d)} \to 0 ~~\text{as}~~ \min_{j \in \{k_1, \ldots, k_m\}} \alpha_j \to \infty ~~\text{while}~~ \max_{j \notin \{k_1, \ldots, k_m\}} \alpha_j = 0,
\end{equation*}
where $k_1, \ldots, k_m$, $m < d$ are the indices where $\beta$ is nonzero.}

Suppose otherwise. Then we find sequences $(\alpha^n)_{n \in \mathbb N} \subset \mathbb N_0^d$ and $(x_n)_{n \in \mathbb N} \subset [0,1]^d$ such that
\begin{equation*}
\min_{j \in \{k_1, \ldots, k_m\}} \alpha^n_j \overset{n \to \infty}{\longrightarrow} \infty ~~\text{and}~~ \forall n \in \mathbb N: \left| \partial^\beta B_\alpha(f)(x_n) - \partial^\beta f(x_n) \right| \ge \epsilon.
\end{equation*}
Since $[0,1]^d$ is compact, the sequence $(x_n)_{n \in \mathbb N}$ has an accumulation point $x_0 \in [0,1]^d$. By the triangle inequality,
\begin{align*}
\left| \partial^\beta B_{\alpha^n}(f)(x_0) - \partial^\beta f(x_0) \right| \ge & \left| \partial^\beta B_\alpha(f)(x_n) - \partial^\beta f(x_n) \right| \\
& - \left| \partial^\beta f(x_0) - \partial^\beta f(x_n) \right| \\
& - \left| \partial^\beta B_{\alpha^n}(f)(x_0) - \partial^\beta B_{\alpha^n}(f)(x_n) \right|.
\end{align*}
\Cref{cor65} and the continuity of $\partial^\beta f$ imply that for any $\delta > 0$, the right hand side is larger than $\epsilon - \delta$ infinitely often. This even contradicts pointwise convergence which is asserted by the inductive hypothesis, and hence the claim is proved.

Now let $\beta \in \mathbb N_0^d$ be arbitrary, and as usual define $\beta' := (0, \beta_2, \ldots, \beta_d)$ and similarly $\alpha' = (0, \alpha_2, \ldots, \alpha_d)$ for each $\alpha \in \mathbb N_0^d$. Then (\Cref{lemma63}, \Cref{thm64})
\begin{align*}
\left\| \partial^\beta B_\alpha(f) - \partial^\beta f \right\| & = \left\| \partial^{\beta'} B_{\alpha'} \left( \partial^{\beta_1 \cdot e_1} B_{\alpha_1 \cdot e_1} (f) \right) - \partial^\beta f \right\| \\
& \le \left\| \partial^{\beta'} B_{\alpha'} \left( \partial^{\beta_1 \cdot e_1} B_{\alpha_1 \cdot e_1} (f) \right) - \partial^{\beta'} B_{\alpha'} \left( \partial^{\beta_1 \cdot e_1} f \right) \right\| + \left\| \partial^{\beta'} B_{\alpha'} \left( \partial^{\beta_1 \cdot e_1} f \right) - \partial^\beta f \right\| \\
& \le \left\| \partial^{\beta_1 \cdot e_1} B_{\alpha_1 \cdot e_1} (f) - \partial^{\beta_1 \cdot e_1} f \right\| + \left\| \partial^{\beta'} B_{\alpha'} \left( \partial^{\beta_1 \cdot e_1} f \right) - \partial^{\beta'} \partial^{\beta_1 \cdot e_1} f \right\|.
\end{align*}
By the claim, both summands on the bottom right go to zero when $\min\limits_{j \in \{1, \ldots, d\}} \alpha_j \to \infty$.
\end{proof}

As a consequence of the chain rule, we have:

\begin{theorem}[Preservation of approximation in derivatives under homothetic transformation]\label{thm68}
For $f, g: \mathbb R^d \to \mathbb R$, assume
\begin{equation*}
\|\partial_\beta f - \partial_\beta g\|_\infty < C.
\end{equation*}
Then
\begin{equation*}
\|\partial_\beta (f \circ h) - \partial_\beta (g \circ h)\|_\infty < \lambda^{|\beta|} C
\end{equation*}
where $h: \mathbb R^d \to \mathbb R^d$ is the homothetic transformation $h(x) = \lambda x + \mathbf v$.
\end{theorem}

Further, we note the following triviality on approximating polynomials by other polynomials:

\begin{theorem}\label{thm69}
Let $p \in \mathbb R[x_1, \ldots, x_d]$, say $p(x) = \sum_{\gamma \le \alpha} a_\gamma x^\gamma$. If $q(x) := \sum_{\gamma \le \alpha} b_\gamma x^\gamma$ is a polynomial such that $|a_\gamma - b_\gamma| < \epsilon$ for all $\gamma \le \alpha$, then
\begin{equation*}
\|\partial_\beta p - \partial_\beta q\|_{\mathcal C([0,1]^d)} \le \epsilon |\alpha| \alpha!.
\end{equation*}
\end{theorem}

We now sharpen \Cref{thm67} a little:

\begin{theorem}
Let $f \in \mathcal C^\gamma(\mathbb R^d)$. Then there exists a sequence $(q_n)_{n \in \mathbb N}$ of rational polynomials such that for $\beta \le \gamma$
\begin{equation*}
\partial^\beta q_n \to \partial^\beta f ~ \text{locally uniformly}, ~~ n \to \infty.
\end{equation*}
\end{theorem}

\begin{proof}
We prove that there exists a sequence of rational polynomials $(q_n)_{n \in \mathbb N}$ such that
\begin{equation*}
\forall n \in \mathbb N: \forall \beta \le \min\{\mathbf n, \gamma\}: \left\|\partial^\beta q_n - \partial^\beta f \right\|_{\mathcal C([-n,n]^d)} < \frac{1}{n}.
\end{equation*}
Indeed, let $n \in \mathbb N$. The homothetic transformation $h(x) = 2nx - \mathbf n$ is a bijection between $[0,1]^d$ and $[-n,n]^d$. Using \Cref{thm67}, pick $p_n \in \mathbb R[x_1, \ldots, x_d]$ such that
\begin{equation*}
\forall \beta \le \min\{\mathbf n, \gamma\}: \left\|\partial^\beta p_n - \partial^\beta (f \circ h) \right\|_{\mathcal C([0,1]^d)} < \frac{1}{2n}.
\end{equation*}
Using \Cref{thm69}, pick $r_n \in \mathbb Q[x_1, \ldots, x_d]$ such that
\begin{equation*}
\forall \beta \le \min\{\mathbf n, \gamma\}: \left\|\partial^\beta r_n - \partial^\beta p_n \right\|_{\mathcal C([0,1]^d)} < \frac{1}{2n}.
\end{equation*}
Now note that by \Cref{thm68},
\begin{equation*}
\left\| \partial^\beta (r_n \circ h^{-1}) - \partial^\beta (f \circ h \circ h^{-1}) \right\|_{\mathcal C([-n,n]^d)} < \frac{1}{n}
\end{equation*}
so that we may pick $q_n := r_n \circ h^{-1}$.
\end{proof}

\bibliographystyle{plain}
\bibliography{multivariate-bernstein}

\begin{thebibliography}{1}

\bibitem{bernstein}
Sergei~Natanovich Bernstein.
\newblock Démonstration du théorème de weierstrass fondée sur le calcul de
  probabilités.
\newblock {\em Communications of the Kharkov Mathematical Society}, 13:1--2,
  1913.

\bibitem{forster}
Otto Forster.
\newblock {\em Analysis 1}.
\newblock Vieweg+Teubner Verlag, 10 edition, 2011.

\bibitem{phillips}
George~M. Phillips.
\newblock {\em Interpolation and Approximation by Polynomials}.
\newblock CMS Books in Mathematics. Springer New York, 2003.

\bibitem{reimer}
Manfred Reimer.
\newblock {\em Multivariate Polynomial Approximation}, volume 144 of {\em
  International Series of Numerical Mathematics}.
\newblock Springer Basel, 2003.

\bibitem{schwartz}
Laurent Schwartz.
\newblock {\em Théorie des distributions}.
\newblock Herrmann, nouv. éd., entièrement corr., refondue et augm. edition,
  1966.

\end{thebibliography}

\end{document}